\documentclass[twoside,10pt,a4paper,leqno]{article}
\usepackage{amssymb,amsmath,amscd,euscript,verbatim,array}
\usepackage[center,pagestyles]{titlesec}
\usepackage{geometry}
\usepackage{anysize}
\usepackage{fancyhdr}
\usepackage{indentfirst}
\usepackage{graphicx}
\usepackage{color}
\usepackage{ifpdf}
\usepackage{mathrsfs}
\usepackage{cite}
\usepackage[numbers,sort&compress]{natbib}
\marginsize{3.5cm}{3.5cm}{2.5cm}{2.5cm}

\titleformat{\section}{\centering\normalsize}{\thesection.}{0.5em}{}
\titleformat{\subsection}{\normalsize\bfseries}{\thesubsection.}{0.5em}{}
\titleformat{\subsubsection}{\normalsize\bfseries}{\thesubsubsection.}{0.5em}{}
\newcommand{\N}{\mathbb{N}}

\newcommand{\R}{\mathbb{R}}

\newtheorem{Theorem}{Theorem}[section]
\newtheorem{Definition}[Theorem]{Definition}
\newtheorem{Lemma}[Theorem]{Lemma}
\newtheorem{Exercise}[Theorem]{Exercise}
\newtheorem{Proposition}[Theorem]{Proposition}
\newtheorem{Remark}[Theorem]{Remark}

\newcommand{\gm}{\gamma}

\newcommand{\T}{\mathbb{T}}
\newcommand{\bthm}{\begin{Theorem}}
\newcommand{\ethm}{\end{Theorem}}
\newcommand{\bpr}{\begin{Proposition}}
\newcommand{\epr}{\end{Proposition}}
\newcommand{\blm}{\begin{Lemma}}
\newcommand{\elm}{\end{Lemma}}
\newcommand{\bex}{\begin{Exercise}}
\newcommand{\eex}{\end{Exercise}}
\newcommand{\be}{\begin{equation}}
\newcommand{\ee}{\end{equation}}
\newcommand{\beal}{\begin{aligned}}
\newcommand{\enal}{\end{aligned}}
\newcommand{\brm}{\begin{Remark}}
\newcommand{\erm}{\end{Remark}}
\newcounter{item}[section]

\newcommand{\Proof}{\textbf{Proof}\hspace{0.3cm}}
\newcommand{\End}{\ensuremath{\hfill{\Box}}\\}
\renewcommand{\title}[1]{\begin{center}\textbf{\large #1}\end{center}}
\renewcommand{\author}[1]{\begin{center}\normalsize #1\end{center}}
\renewcommand{\date}[1]{\begin{center}#1\end{center}}

\makeatletter \@addtoreset{equation}{section}
\makeatother
\begin{document}
\vspace{10pt}
\title{ASYMPTOTIC LIPSCHITZ REGULARITY OF VISCOSITY SOLUTIONS OF HAMILTON-JACOBI EQUATIONS   }

\vspace{6pt}
\author{\sc Xia Li and Lin Wang}

\vspace{10pt} \thispagestyle{plain}

\begin{quote}
\small {\sc Abstract.} For each continuous initial data $\varphi(x)\in C(M,\R)$, we obtain the asymptotic Lipschitz regularity  of the viscosity solution of the following evolutionary Hamilton-Jacobi equation with convex and coercive Hamiltonians:
\begin{equation*}
\begin{cases}
\partial_tu(x,t)+H(x,\partial_xu(x,t))=0,\\
u(x,0)=\varphi(x).
\end{cases}
\end{equation*}
\end{quote}
\begin{quote}
\end{quote}
\begin{quote}
\end{quote} \vspace{25pt}

\section{\sc Introduction and main result}

Let $M$ be a $n$-dimensional connected and closed smooth manifold. We are concerned with a  Hamiltonian $H:T^*M\rightarrow \R$ satisfying the following assumptions:
\begin{itemize}
\item [\textbf{(H1)}] \textbf{Smoothness}: $H(x,p)$ is a $C^2$ function;
\item [\textbf{(H2)}] \textbf{Convexity}: $H(x,p)$ is strictly convex with respect to $p$;
\item [\textbf{(H3)}] \textbf{Coercivity}: for each  $x\in M$, $H(x,p)\rightarrow\infty$  uniformly as $|p|\rightarrow\infty$.
\end{itemize}
(H3) is equivalent to the topological statement that for  each $c\in\R$, the set $\{(x,p)\in T^*M|x\in K,H(x,p)\leq c\}$ is compact.

We consider the following  Hamilton-Jacobi equation under the assumptions (H1)-(H3):
\begin{equation}\label{hje1}
\begin{cases}
\partial_tu(x,t)+H(x,\partial_xu(x,t))=0,\\
u(x,0)=\varphi(x),
\end{cases}
\end{equation}where $(x,t)\in M\times[0,\infty)$ and $\varphi(x)\in C(M,\R)$.

We recall  the Ma\~{n}\'{e} critical value of $H(x,p)$ denoted by $c[0]$. By \cite{CIPP}, one has
\begin{equation}
c[0]=\inf_{u\in C^1(M,\R)}\max_{x\in M}H(x,\partial_xu).
\end{equation}
Let $u(x,t)$ be the viscosity solution of (\ref{hje1}). It was shown by \cite{ds} that the the limit $v(x):=\lim_{t\rightarrow \infty}(u(x,t)+c[0]t)$ is a Lipschitz weak KAM solution of
\begin{equation}\label{hje}
H(x,\partial_xu)=c[0].
\end{equation}
Recently, a convergence result for more general contact Hamilton-Jacobi equations was established in \cite{SWY}. Note that the limit $v(x)$ is a Lipschitz function, while the initial data $\varphi(x)$ is only continuous. A  question is

{\it when does the Lipschitz regularity of the viscosity solution of (\ref{hje1}) emerge?}

If $H(x,p)$ is superlinear with respect to $p$, then the Lipschitz regularity emerges after an arbitrarily small time, which is basically from the celebrated Fleming's lemma \cite[Theorem 4.4.3]{F3}. Unfortunately, if $H(x,p)$ is coercive, the Fleming's lemma does not hold anymore. Then it is natural to ask that

{\it will the Lipschitz regularity of the viscosity solution of (\ref{hje1}) emerge after a finite time (asymptotic Lipschitz regularity) or  an infinite time (limit Lipschitz regularity)?}

In this note, we clarify the asymptotic Lipschitz regularity of the viscosity solution of (\ref{hje1}) is true.
 More precisely, we have:
\begin{Theorem}\label{two}
Let $u(x,t)$ be a viscosity solution of (\ref{hje1}) with continuous initial data $\varphi\in C(M,\R)$, then  there exists $t_0>0$ such that for $t>t_0$, $u(x,t)$ is  $\iota$-Lipschitz continuous, where $t_0, \iota:=\iota(t_0)$ are independent of $\varphi$.
\end{Theorem}

This note is outlined as follows. In Section 2, some properties of viscosity solutions are introduced as preliminaries. In Section 3, by introducing a modified Hamiltonian, the Ma\~{n}\'{e} critical value and action minimizing orbits are located. The proof of Theorem \ref{two} is completed in Section 4.

\section{\sc Preliminaries}
In this section,   we introduce some properties of the viscosity solutions in our settings.
First of all, we introduce the notion of  semiconcave functions.

\begin{Definition}[Semiconcavity on $\R^n$]Let $U$ be an open convex subset of $\R^n$ and let $u:U\rightarrow\R$ be a function. $u$ is called a semiconcave function with linear modulus if there exists a finite constant $K$ and for each $x\in U$ there exists a linear form $\theta_x:\R^n\rightarrow\R$ such that for any $y\in U$
\begin{equation}
u(y)-u(x)\leq \theta_x(y-x)+K|y-x|^2.
\end{equation}
\end{Definition}
For the sake of simplicity, we only consider the semiconcave functions with linear modulus defined as above. See \cite{CS} for a more general definition. In this context, the notion ``semiconcave" means ``semiconcave with a linear modulus".
\begin{Definition}[Semiconcavity on a manifold]
A function $u:M\rightarrow\R$ defined on the $C^r$ ($r\geq 2$) differential $k$-dimensional manifold $M$ is locally semiconcave if for each $x\in M$ there exists a $C^r$ ($r\geq 2$) coordinate chart $\psi:U\rightarrow\R^n$ with $x\in U$ such that $u\circ \psi^{-1}:U\rightarrow\R$ is semiconcave.
\end{Definition}
Consider the stationary equation
\begin{equation}\label{hstaee}
H(x,\partial_xu)=0,
\end{equation}
and the evolutionary equation
\begin{equation}\label{hevaee}
\partial_tu+H(x,\partial_xu)=0,
\end{equation}
based on \cite[Theorem 5.3.1. and Theorem 5.3.6]{CS}, we have the following results.
\begin{Proposition}\label{semi}
Let $H\in C^2(T^*M,\R)$, we have the following properties.
\begin{itemize}
\item [(a)] Let $u$ be a semiconcave function satisfying the equations (\ref{hstaee}) (resp. (\ref{hevaee})) almost everywhere. If $H(x,p)$ is  convex with respect to $p$, then $u$ is a viscosity solution of the equations (\ref{hstaee}) (resp. (\ref{hevaee}));
\item [(b)] Let $u$ be a Lipschitz viscosity solution of the equations (\ref{hstaee}) (resp. (\ref{hevaee})). If $H(x,p)$ is strictly convex with respect to $p$, then $u$ is locally semiconcave on $M$ (resp. $M\times(0,+\infty)$).
\end{itemize}
\end{Proposition}

Let us recall the notion of upper differentials (see \cite{CS,F3} for instance).
\begin{Definition}[Upper differential on $\R^n$]
Let $u:U\rightarrow\R$ be a function defined on the open subset $U$ of $\R^n$, the set
\[D^+u(x_0):=\left\{\theta\in \R^n\ \big|\ \limsup_{x\rightarrow x_0}\frac{u(x)-u(x_0)-\theta(x-x_0)}{|x-x_0|}\leq 0\right\}\]
 is called a upper  differential of $u$ at $x_0$.
\end{Definition}
\begin{Definition}[Upper differential on a manifold]
Let $u:M\rightarrow\R$ be a function defined on the the manifold $M$, the linear form $\theta\in T^*_{x_0}M$ is a upper differential of $u$ at $x_0\in M$, if there exist a neighborhood $V$ of $x_0$ and a function $\varphi:V\rightarrow\R$, diffferentiable at $x_0$, with $\varphi(x_0)=u(x_0)$ and $d_{x_0}\varphi=\theta$ and such that $\varphi(x)\geq u(x)$ for each $x\in V$.
\end{Definition}
It is easy to verify the equivalence between the definition of upper differentials on an Euclidean space and the one on a manifold.

We use $\partial u(x_0,\theta)$ to denote one-sided directional derivative along $\theta\in \R^n$ at $x_0$, namely
\[\partial u(x_0,\theta):=\lim_{h\rightarrow 0^+}\frac{u(x_0+h\theta)-u(x_0)}{h}.\]
The upper differential and  one-sided directional derivative of the semiconcave function enjoy the following properties (\cite[Proposition 3.3.4 and Theorem 3.3.6]{CS}).
\begin{Proposition}\label{semiapp}
Let $u:M\rightarrow\R$ be a semiconcave function. Then following properties hold true.
\begin{itemize}
\item [(a)] $D^+u(x)\neq \emptyset$ for any $x\in M$;
\item [(b)] If $\{x_n\}$ is a sequence in $M$ converging to $x$ and if $p_n\in D^+u(x_n)$ converges to a vector $p$, then $p\in D^+u(x)$;
\item [(c)] $\partial u(x,\theta)=\min_{p\in D^+u(x)}\langle p,\theta\rangle$ for any $x\in M$ and $\theta\in \R^n$.
\end{itemize}
\end{Proposition}

Throughout this paper, we shall use $| \cdot |$ to denote the Euclidean norm, that is  $|\alpha|=\sqrt{\alpha^2_1+\ldots+\alpha^2_i}$ for given $\alpha=(\alpha_1,\ldots,\alpha_i)\in \R^i$, $i=1$ or $i=n$.
\section{\sc Ma\~{n}\'{e} critical value and action minimizing orbits}

\subsection{Modification of the Hamiltonian}

Let $H(x,p)$ be a Hamiltonian satisfying (H1)-(H3). We construct a new Hamiltonian denoted by $H_{R}(x,p)$ with $R>1$ as follows. Without loss of generality, we assume $M=\T^n$, from which $T^*M=\T^n\times\R^n$.
\begin{equation}
H_{R}(x,p)=\alpha_R (p)H(x,p)+\mu_R\beta(|p|^2-R^2),
\end{equation} where $\mu_R$ is a  constant determined by (\ref{murr}) below and $\alpha_R (p)$ is a $C^2$ function satisfying
\begin{equation}
\alpha_R (p)=\left\{\begin{array}{ll}
\hspace{-0.4em}1,&  |p|\leq R+1,\\
\hspace{-0.4em}0,&|p|>R+2.\\
\end{array}\right.
\end{equation}
Without loss of generality, one can require $|\alpha'_R (p)|<2$ and $\|\alpha''_R (p)\|_1<2$, where $\|\cdot\|_1$ denotes 1-norm, namely the maximum of the summation of the absolute values of elements in each column.
$\beta(z) $ is defined as \begin{equation}
\beta(z)=\left\{\begin{array}{ll}
\hspace{-0.4em}0,&  |z|\leq 0,\\
\hspace{-0.4em}z^4,&|z|>0,\\
\end{array}\right.
\end{equation}
It is easy to see that $H_R(x,p)=H(x,p)$ for $|p|\leq R$.
In the following, we show that $H_R(x,p)$  satisfies (H1),(H2) and superlinearity.

\vspace{1em}

\noindent\textbf{Claim 1.} {\it $H_R(x,p)$  satisfies (H1).}

\noindent\textbf{Proof of Claim 1.}
Note that $\alpha_R(p)$ and $H(x,p)$ are $C^2$ functions. By the construction, $\beta(z)$ is of class $C^3$. It follows that $H_R(x,p)$ is a $C^2$ function.
\End

\vspace{1em}

\noindent\textbf{Claim 2.} {\it $H_{R}(x,p)$  satisfies (H2).}

\noindent\textbf{Proof of Claim 2.}
It suffices to show that for given $x\in M$,  $\partial^2 H_R/\partial {p}^2 (x,p)>0$.
\begin{itemize}
\item [(i)] For $|p|\leq R$,
\[H_R(x,p)=H(x,p).\]
Hence, we have
\[\frac{\partial^2 H_R}{\partial {p}^2} (x,p)=\frac{\partial^2 H}{\partial {p}^2} (x,p)>0.\]
\item [(ii)] For $R<|p|\leq R+1$,
\begin{align*}
H_R(x,p)=H(x,p)+\mu_R\beta(|p|^2-R^2).
\end{align*}
It follows that
\[\frac{\partial^2 H_R}{\partial {p}^2} (x,p)=\frac{\partial^2 H}{\partial {p}^2} (x,p)+2\mu_R\left(2\beta''(|p|^2-R^2) Z(p)+\beta'(|p|^2-R^2)\cdot E\right)>0,\]
where $Z(p):=(p_1,\ldots,p_n)^T\cdot (p_1,\ldots,p_n)$, $E$ denotes the $n\times n$ identity matrix.
\item [(iii)] For $R+1<|p|\leq R+2$,
\begin{align*}
H_R(x,p)=\alpha_R (p) H(x,p)+\mu_R\beta(|p|^2-R^2).
\end{align*}
It yields that
\begin{align*}
\frac{\partial^2 H_R}{\partial {p}^2} (x,p)=&H(x,p)\alpha''_R (p)+W(x,p)+
\alpha_R (p)\frac{\partial^2 H}{\partial {p}^2} (x,p)\\
&+2\mu_R\left(2\beta''(|p|^2-R^2)Z(p)+\beta'(|p|^2-R^2)\cdot E\right),
\end{align*}
where \[W(x,p):=\alpha'_R (p)^T\cdot\frac{\partial H}{\partial {p}} (x,p)+\frac{\partial H}{\partial {p}} (x,p)^T\cdot\alpha'_R (p).\]
Since $W(x,p)$ is symmetric, then $\partial^2 H_R/\partial {p}^2 (x,p)$ is symmetric. We denote $\partial^2 H_R/\partial {p}^2 (x,p)=(a_{ij})_{n\times n}$, then $\partial^2 H_R/\partial {p}^2 (x,p)(x,p)$ is positive definite if $\sqrt{a_{ii}a_{jj}}>(n-1)|a_{ij}|$ for $i,j=1,\ldots,n$ and $i\neq j$.

Based on the construction of $\alpha_R (p)$ and the compactness of $M$, let
\[\gamma_R:=2\max_{(x,p)\in \T^n\times[R+1,R+2]^n}|H(x,p)|+(n-1)\max_{(x,p)\in \T^n\times[R+1,R+2]^n}\|W(x,p)\|_1,\]
it is enough to take
\begin{equation}\label{murr}
\mu_R>\max\left\{\gamma_R
 , 1\right\}.
 \end{equation}
\item [(iv)] For $|p|> R+2$,
\[H_R(x,p)=\mu_R\beta(|p|^2-R^2),\]
which implies
\[\frac{\partial^2 H_R}{\partial {p}^2} (x,p)=2\mu_R\left(2\beta''(|p|^2-R^2) Z(p)+\beta'(|p|^2-R^2)\cdot E\right)>0\]
\end{itemize}
Therefore, $H_{R}(x,p)$  satisfies (H2).
\End

\vspace{1em}

\noindent\textbf{Claim3.} {\it $H_{R}(x,p)$  satisfies the superlinearity.}

\noindent\textbf{Proof of Claim 3.} It suffices to verify the superlinearity of $H_R(x,p)$ for $|p|>R+2$. In this case, we have
\[H_R(x,p)\geq \mu_R\beta(|p|^2-R^2)\geq |p|^2.\]
Hence, for each $A>0$, one can find $C_A>0$ such that
\[H_R(x,p)\geq A|p|-C_A.\]
Therefore, $H_{R}(x,p)$  satisfies the superlinearity.
\End

 It is easy to see that
$H_R$ converges uniformly on compact subsets to $H$ in the $C^2$ topology as $R\rightarrow\infty$.

\subsection{Ma\~{n}\'{e} critical value}

We use $c_R$ to denote the Ma\~{n}\'{e} critical value of $H_R(x,p)$. Then
\begin{equation}\label{criti}
c_R=\inf_{u\in C^1(M,\R)}\max_{x\in M}H_R(x,\partial_xu).
\end{equation}
The following lemma asserts that for $R$ large enough, the Ma\~{n}\'{e} critical value of $H_R$ is independent of $R$. We  denote
\begin{equation}
c[0]:=\inf_{u\in C^1(M,\R)}\max_{x\in M}H(x,\partial_xu),
\end{equation}
which can be seen as the Ma\~{n}\'{e} critical value of $H(x,p)$.
\begin{Lemma}\label{cr=cr}
There exists $R_0>0$ such that for any $R>R_0$, we have
\begin{equation}\label{crr}
c_R=c[0].
\end{equation}
\end{Lemma}

\noindent\Proof
From (\ref{criti}) and the construction of $H_R$, it follows that for any $R>0$,
\begin{equation}
c_R\leq \max_{x\in M}H_R(x,0)=\max_{x\in M}H(x,0).
\end{equation}
Let $A:=\max_{x\in M}H(x,0)+1$. We denote
\[\Lambda:=\{(x,p)\in T^*M|x\in M,H(x,p)\leq A\}.\]
By (H3) and the compactness of $M$, $\Lambda$ is compact. Hence, there exists $R_0>0$ such that
 \[\Lambda\subset \{(x,p)\in T^*M|x\in M,|p|_x\leq R_0\},\]
 where $|\cdot|_x$ denotes the Riemannian metric on $T_x^*M$. Based on the construction of $H_R$, it yields that for any $R>R_0$ and $(x,p)\in \Lambda$, we have
\begin{equation}
H_R(x,p)=H(x,p).
\end{equation}
In terms of the definition of the Ma\~{n}\'{e} critical value, one can find a sequence $u_n\in C^1(M,\R)$ such that
\begin{equation}\label{cr}
\max_{x\in M}H_R(x,\partial_xu_n(x))\rightarrow c_R.
\end{equation}Since $c_R< A$, then we have $|\partial_xu_n(x)|\leq R_0$ for $n$ large enough. Moreover, we have $H_R(x,\partial_xu_n(x))=H(x,\partial_xu_n(x))$ for any $R>R_0$. Then, it yields for $n$ large enough,
\begin{align*}
c[0]&=\inf_{u\in C^1(M,\R)}\max_{x\in M}H(x,\partial_xu(x)),\\
&\leq \max_{x\in M}H(x,\partial_xu_n(x)),\\
&= \max_{x\in M}H_R(x,\partial_xu_n(x)).
\end{align*}
Taking the limit as $n\rightarrow\infty$, it follows from (\ref{cr}) that $c[0]\leq c_R$. Similarly, we choose
a sequence $v_n\in C^1(M,\R)$ such that
\begin{equation}\label{cr1}
\max_{x\in M}H(x,\partial_xv_n(x))\rightarrow c[0].
\end{equation}Since $c[0]\leq \max_{x\in M} H(x,0)<A$, then we have $|\partial_xv_n(x)|\leq R_0$ for $n$ large enough. Moreover, we have $H_R(x,\partial_xv_n(x))=H(x,\partial_xv_n(x))$ for any $R>R_0$. Then, it yields for $n$ large enough,
\begin{align*}
c_R&=\inf_{u\in C^1(M,\R)}\max_{x\in M}H_R(x,\partial_xu(x)),\\
&\leq \max_{x\in M}H_R(x,\partial_xv_n(x)),\\
&= \max_{x\in M}H(x,\partial_xv_n(x)),
\end{align*}
which together with (\ref{cr1}) implies that $c_R\leq c[0]$ as $n\rightarrow\infty$.
 Therefore, one can find $R_0>0$ such that for any $R>R_0$,
$c_R=c[0]$. This finishes the proof of Lemma \ref{cr=cr}.
\End

For the sake of simplicity,  we assume $c[0]=0$ in the following context.
\subsection{The viscosity solution of (\ref{hje})}
Let $\bar{u}(x)$ be a viscosity solution of $H(x,\partial_xu)=0$. Since $H(x,p)$ is coercive with respect to $p$, then  $\bar{u}(x)$ is a  Lipschitz function on $M$, which together with Proposition \ref{semi} implies that $\bar{u}$ is  semiconcave.

Let $\mathcal{D}$ be the set of all differentiable points of $\bar{u}$ on $M$. Due to the Lipschitz property of $\bar{u}$, it follows that $\mathcal{D}$ has full Lebesgue measure.

\begin{Lemma}\label{th}
There exists $R_1>0$ such that for any $R>R_1$, $\bar{u}(x)$ is a viscosity solution of $H_R(x,\partial_xu)=0$.
\end{Lemma}
\Proof
Since  $\bar{u}(x)$ is a  Lipschitz function on $M$, then for $x\in \mathcal{D}$, we have $H(x,\partial_x\bar{u})=0$. By (H3), there exists $R_1>0$ such that $|\partial_x\bar{u}|\leq R_1$ for $x\in \mathcal{D}$. It follows from the construction of $H_R$ that for $R>R_1$ and
 \[(x,p)\in \{(x,p)\in T^*M|x\in \mathcal{D},|p|_x\leq R_1\},\]
 we have $H_R(x,p)=H(x,p)$, which means that for $x\in \mathcal{D}$,
\[H_R(x,\partial_x\bar{u})=0.\]
Due to the semiconcavity of $\bar{u}(x)$, it follows from Proposition  \ref{semi}  that $\bar{u}(x)$ is a viscosity solution of $H_R(x,\partial_xu)=0$ for any $R>R_1$. This completes the proof of Lemma \ref{th}.
\End

\subsection{Location of the action minimizing orbits}
Let $\Phi^t_H$ denote the flow  generated by $H(x,p)$. Let $(x(t),p(t)):=\Phi^t_H(x_0,p_0)$.
Let $L_R$ be the Lagrangian associated to $H_R$. To fix the notion, for a given $R>0$ and $(x_0,p_0)\in T^*M$, we call $(x_R(t),p_R(t)):=\Phi_{H_R}^t(x_0,p_0)$ the action minimizing orbit with  $x_R(0)=x_0$ and $x_R(t)=y$ if
\[x_R(t)=\gm_R(t),\quad p_R(t)=\frac{\partial L_R}{\partial \dot{x}}(\gm_R(t),\dot{\gm}_R(t)),\]
where $\gm_R:[0,t]\rightarrow M$ is an action minimizing curve with  $\gm_R(0)=x_0$ and $\gm_R(t)=y$. That is $\gm_R$ achieves
\[\inf_{\substack{\gm(0)=x_0\\ \gm(t)=y}}\int_0^tL_R(\gm(s),\dot{\gm}(s))ds.\]

\begin{Lemma}[{\it a priori} compactness]\label{key}
For $s\in [0,t]$, let $(x_R(s),p_R(s))$ be an action minimizing orbit with $x_R(0)=x_0$ and $x_R(t)=y$.  There exists $\bar{R}>1$ such that for any $R>\bar{R}$, one can find $t_0:=t_0(\bar{R})>0$ such that for any  $s\in [0,t]$ with $t>t_0$, we have
\[(x_R(s),p_R(s))\in \Omega,\]
where $\Omega:=\{(x,p)\ |\ H(x,p)\leq 1\}$.
\end{Lemma}

In order to prove Lemma \ref{key}, we need to do some preparations.
Based on Lemma \ref{th}, it yields that  for $x\in \mathcal{D}$ and $R>R_1$,
\begin{equation}\label{hxud}
H_R(x, \partial_x \bar{u}(x) ) = 0.
\end{equation}
We define
\begin{equation}\label{define}
\widetilde{L}_R(x,\dot{x})= L_R(x, \dot{x}) - \langle \partial_x\bar{u}(x), \dot{x}\rangle,\quad  x\in \mathcal{D}.
\end{equation}
Denote
\begin{equation}\label{gamm}
\Gamma_R :=\left\{ \left(x, \frac{\partial H_R}{\partial p}(x,\partial_x \bar{u}(x))\right)~:~x\in \mathcal{D}  \right\},
\end{equation}
where $\frac{\partial H_R}{\partial p}$ denotes the partial derivative of $H_R$ with respect to the second argument. We have the following lemma.
\begin{Lemma}\label{squr}
For any $x\in\mathcal{D}$, $\widetilde{L}_R(x,\dot{x})\geq 0$. In particular, $\widetilde{L}_R(x,\dot{x})=0$ if and only if $(x,\dot{x})\in \Gamma_R$.
\end{Lemma}
\Proof
By (\ref{define}) and (\ref{gamm}), we have
\begin{equation}\label{11}
\widetilde{L}_R\bigg|_{\Gamma_R} =-H_R(x, \partial_x \bar{u}(x))= 0.
\end{equation}
In addition, we have
\begin{equation}\label{12}
\frac{\partial \widetilde{L}_R}{\partial \dot{x}}\bigg|_{\Gamma_R}
= \frac{\partial L_R}{\partial \dot{x}}(x,\dot{x}) - \partial_x \bar{u}(x) = 0.
\end{equation}

By the superlinearity of $L_R$, it follows  from (\ref{11}) that there exists $K_1>0$ large enough such that for $|\dot{x}|> K_1$,
\[\widetilde{L}_R (x, \dot{x})\geq d>0,\]where $d$ is a constant independent of $(x,\dot{x})$.

For $x\in \mathcal{D}$, $\bar{u}(x)$ satisfies the equation (\ref{hxud}). Since $\bar{u}(x)$ is Lipschitz continuous, then  $\partial_x\bar{u}(x)$ is bounded. Let
\[\dot{x}_0:=\frac{\partial H_R}{\partial p}(x,\partial_x \bar{u}(x)),\]
 then there exists $K_2>0$ independent of $x$ such that $|\dot{x}_0|\leq K_2$. Take $K_3:=\max\{K_1,K_2\}$. Note that
$\frac{\partial^2L_R}{\partial \dot{x}^2}(x,\dot{x})$ is positive definite, for $|\dot{x}|\leq K_3$, it follows  from (\ref{11}) and (\ref{12}) that there exists $\Lambda>0$ independent of $(x,\dot{x})$ such that
\begin{equation}
\widetilde{L}_R (x, \dot{x}) \geq \Lambda\left|\dot{x} - \frac{\partial H_R}{\partial p}(x,\partial_x \bar{u}(x))\right|^2.
\end{equation}
Consequently, it is easy to see that
\begin{equation}\label{Lagrangian graph}
\widetilde{L}_R(x, \dot{x})
\left\{\!\!\!
  \begin{array}{rl}
   &=0,  \qquad (x,\dot{x})\in \Gamma_R,\\
   &>0,  \qquad (x,\dot{x}) \notin \Gamma_R.
  \end{array}
\right.
\end{equation}This completes the proof of Lemma \ref{squr}.
\End

Let $\Omega^*$ denote the Legendre transformation of $\Omega$ via $\mathcal{L}: T^*M\rightarrow TM$. By (H3), there exist $R_2, R^*_2>0$ such that
\[\Omega\subset\{(x,p)\in T^*M\ |\ x\in M,\ |p|_x\leq R_2\}.\]
\[\Omega^*\subset\{(x,v)\in TM\ |\ x\in M,\ |v|_x\leq R^*_2\}.\]
Based on the preparations above, we will prove Lemma \ref{key}. First of all, we take
\begin{equation}\label{rrr}
\bar{R}=\max\{R_0,R_1,R_2,R^*_2\},
\end{equation} where $R_0$, $R_1$ are determined by Lemma \ref{cr=cr} and Lemma \ref{th}.

\textbf{Proof of Lemma \ref{key}:}
By the energy conservation of $H$, it suffices to prove $(x_0,p_0)\in \Omega$, where $(x_0,p_0)=(x_R(0),p_R(0))$ is the initial point of the flow $\Phi_{H_R}^t$.
Let
\begin{equation}
\Delta:=T^*M\backslash \Omega=\{(x,p)\ |\ H(x,p)> 1\}.
\end{equation}By contradiction, we assume $(x_0,p_0)\in \Delta$.

Let $\Sigma:=\{(x, \partial_x \bar{u}(x))\ |\ x\in \mathcal{D}\}$.
Since $H(x, \partial_x \bar{u}(x))= 0$ for $x\in \mathcal{D}$, then $\Sigma\cap \Delta=\emptyset$. Let $\Sigma^*$ and $\Delta^*$ denote the Legendre transformation of $\Sigma$ and $\Delta$ via $\mathcal{L}: T^*M\rightarrow TM$ respectively. Since $\mathcal{L}$ is a diffeomorphism onto the image, then we have
\begin{equation}\label{222}
\Sigma^*\cap \Delta^*=\emptyset.
\end{equation} By virtue of Lemma \ref{th}, it yields that for $R>\bar{R}$ and $ x\in \mathcal{D}$,
\[\frac{\partial H_R}{\partial p}(x,\partial_x \bar{u}(x))=\frac{\partial H}{\partial p}(x,\partial_x \bar{u}(x)).\] It follows that
\begin{equation}
\Sigma^*=\left\{\left(x, \frac{\partial H}{\partial p}(x,\partial_x \bar{u}(x))\right)\ :\ x\in \mathcal{D}\right\}.
\end{equation}

We use $\Sigma^*_\kappa$ to denote a $\kappa$-neighborhood of $\Sigma^*$ in the fibers, namely
\[\Sigma^*_\kappa:=\left\{(x,\dot{x})\ \big|\ x\in \mathcal{D}, \text{dist}\left(\dot{x}, \frac{\partial H}{\partial p}(x,\partial_x \bar{u}(x))\right)\leq \kappa\right\}.\]

By the $C^2$ regularity of $H$ and $L_R$,  for any $\epsilon>0$, there exists $\kappa>0$ such that for $(x,\dot{x})\in \Sigma^*_\kappa$, we have
 \[H\left(x,\frac{\partial L_R}{\partial \dot{x}}(x,\dot{x})\right)\leq \epsilon,\]
hence, for $\kappa$ small enough, we have $\epsilon<1$. Moreover
 \[\Sigma^*_\kappa\cap\Delta^*=\emptyset.\]
By Lemma \ref{squr}, for any $x\in\mathcal{D}$, $\widetilde{L}_R(x,\dot{x})\geq 0$ and $\widetilde{L}_R(x,\dot{x})=0$ if and only if $(x,\dot{x})\in \Sigma^*$. Then  for each $R>\bar{R}$ , there exists  a constant $\eta:=\eta(\bar{R})>0$ such that for $x\in \mathcal{D}$ and $(x,\dot{x})\in \Delta^*$
\begin{equation}\label{eta1}
\widetilde{L}_R(x, \dot{x})\geq \eta,
\end{equation}
where
\[\widetilde{L}_R(x,\dot{x})= L_R(x, \dot{x}) - \langle \partial_x\bar{u}(x), \dot{x}\rangle.\]

Let $\gm_R:[0,t]\rightarrow M$ be an action minimizing curve with  $\gm_R(0)=x_0$, $\gm_R(t)=y$. Then we have $\dot{\gm}_R(s)=\frac{\partial H_R}{\partial p}(x_R(s),p_R(s))$ for $s\in [0,t]$. Since  $(x_0,p_0)\in \Delta$, then for $s\in [0,t]$, we have
\begin{equation}\label{eta4}
(\gm_R(s),\dot{\gm}_R(s))\in \Delta^*.
\end{equation}
Let $\Theta$ be the set of $\gm_R(s)$ along which the one-sided directional derivative denoted by $\partial \bar{u}(\gm_R(s), \dot{\gm}_R(s))$ exists. For $\gm_R(s)\in \Theta$, we denote
\[\widehat{L}_R(\gm_R(s), \dot{\gm}_R(s)):=L_R(\gm_R(s), \dot{\gm}_R(s))-\partial \bar{u}(\gm_R(s), \dot{\gm}_R(s)).\]
Note that $\bar{u}$  is locally semiconcave. By virtue of Proposition \ref{semiapp} (b), one can find a sequence $x_n^s\in \mathcal{D}$ with $x_n^s\rightarrow \gm_R(s)$ and $\partial_x\bar{u}(x_n^s)\rightarrow p_s\in D^+\bar{u}(\gamma_R(s))$ as $n\rightarrow\infty$ for a given $s\in [0,t]$. By virtue of Proposition \ref{semiapp} (c),  for $n$ large enough, extracting a subsequence if necessary, we have
\begin{align*}
\partial \bar{u}(\gm_R(s), \dot{\gm}_R(s))&=\min_{p\in D^+\bar{u}(\gamma_R(s))}\langle p,\dot{\gamma}_R(s)\rangle,\\
&\leq \langle p_s,\dot{\gamma}_R(s)\rangle,\\
&\leq \langle \partial_x\bar{u}(x_n^s),\dot{\gm}_R(s)\rangle+\frac{1}{n}.
\end{align*}
Note that $\Delta^*$ is an open set, then $(x_n^s,\dot{\gamma}_R(s))\in \Delta^*$ for $n$ large enough. It follows from  (\ref{eta1}) that  for every $s\in [0,t]$ and $n$ large enough,
\begin{equation}\label{eta}
\widehat{L}_R(\gm_R(s), \dot{\gm}_R(s))\geq L_R(x_n^s,\dot{\gm}_R(s))-\langle \partial_x\bar{u}(x_n^s),\dot{\gm}_R(s)\rangle-\frac{2}{n}\geq \frac{\eta}{2}.
\end{equation}
Moreover, we have
\begin{align*}
\int_0^t\widehat{L}_R(\gm_R(s), \dot{\gm}_R(s))ds\geq \frac{\eta}{2}t.
\end{align*}
On the other hand, we have
\begin{align*}
\int_0^t\widehat{L}_R(\gm_R(s), \dot{\gm}_R(s))ds&=\int_0^tL_R(\gm_R(s), \dot{\gm}_R(s))-\partial \bar{u}(\gm_R(s), \dot{\gm}_R(s))ds,\\
&=\int_0^tL_R(\gm_R(s), \dot{\gm}_R(s))ds-(\bar{u}(\gm_R(t))-\bar{u}(\gm_R(0))).
\end{align*}
It follows from the semiconcavity and  the compactness of $M$ that $\bar{u}$ has a uniform bound denoted by $C_0$. Hence, we have
\begin{equation}\label{eta2}
\int_0^tL_R(\gm_R(s), \dot{\gm}_R(s))ds\geq \frac{\eta}{2}t-2C_0.
\end{equation}

On the other hand,  $\gm_R$ is an action minimizing curve of $L_R$. Let $\gm_{R_2}$ be an action minimizing curve of $L_{R_2}$. It follows from Lemma \ref{cr=cr}   that for $R>\bar{R}$, there exists a constant $C_1>0$ independent of $R$ such that
\begin{equation}\label{ccrr}
\begin{split}
\int_0^tL_R(\gm_R(s), \dot{\gm}_R(s))ds&\leq \int_0^tL_{R}(\gm_{R^*_2}(s), \dot{\gm}_{R^*_2}(s))ds,\\
&=\int_0^tL_{R^*_2}(\gm_{R^*_2}(s), \dot{\gm}_{R^*_2}(s))ds,\\
&=h^t_{R^*_2}(x_0,y)\leq C_1,
\end{split}
\end{equation}
 where $h^t_{R^*_2}(x_0,y)$ denotes the minimal action of $L_{R^*_2}$. It is clear to see that (\ref{ccrr}) contradicts with (\ref{eta2}) if we take $t>(4C_0+2C_1)/\eta$. Let $t_0:=(4C_0+2C_1)/\eta$, then we have  $(x_0,p_0)\notin \Delta$ for $t>t_0$. Obviously, $t_0$ only depends on $\bar{R}$. This completes the proof of Lemma \ref{key}.
\End

\section{\sc Asymptotic Lipschitz regularity}
In this section, we are devoted to proving Theorem \ref{two}, which is concerned with  the following  Hamilton-Jacobi equation under the assumptions (H1)-(H3):
\begin{equation}\label{hje11}
\begin{cases}
\partial_tu(x,t)+H(x,\partial_xu(x,t))=0,\\
u(x,0)=\varphi(x),
\end{cases}
\end{equation}where $(x,t)\in M\times[0,\infty)$ and $\varphi(x)\in C(M,\R)$. Let $u_{R}(x,t)$ be the viscosity solution of the following equation:
\begin{equation}\label{hra}
\begin{cases}
\partial_tu(x,t)+H_{R}(x,\partial_xu(x,t))=0,\\
u(x,0)=\varphi(x).
\end{cases}
\end{equation}
Let $T_t^{R}$ be the Lax-Oleinik semigroup generated by $L_R$ associated to $H_R$ via the Legendre transformation. Namely,
\begin{equation}\label{semmmmm}
T_t^{R}\varphi(x)=\inf_{\gm(t)=x}\left\{\varphi(\gm(0))+\int_0^tL_R(\gm(s),\dot{\gm}(s))ds\right\}.
\end{equation}
 Then we have
\begin{equation}\label{lo}
u_R(x,t)=T_t^{R}\varphi(x).
\end{equation}
First of all, we consider the viscosity solutions of (\ref{hje11}) with $t$ suitable large.
\begin{Lemma}\label{tgeq}
For any $R\geq \bar{R}$ where $\bar{R}$ is determined by (\ref{rrr}), there exists $t_0>0$ such that for $t>t_0$, $u_{R}(x,t)$ is a viscosity solution of the following equation:
\begin{equation}\label{uh}
\partial_tu(x,t)+H(x,\partial_xu(x,t))=0.
\end{equation}
\end{Lemma}
\Proof
By Proposition \ref{semi} (b), $u_{R}(x,t)$ is locally semiconcave on $M\times(0,\infty)$. Let $\mathcal{E}_R$ be the set of all differentiable points of $u_{R}(x,t)$ on $M\times(0,\infty)$, then $\mathcal{E}_R$ has full Lebesgue measure. For $(x,t)\in \mathcal{E}_R$, we have $u_{R}(x,t)$ satisfies (\ref{hra}). For a given $(\bar{x},\bar{t})\in \mathcal{E}_R$, let $\gm_{R}:[0,\bar{t}]\rightarrow M$ be a curve  achieving the infimum of (\ref{semmmmm}) with $\gm_{R}(\bar{t})=\bar{x}$. Then we have
\begin{equation}
\partial_xu_{R}(\bar{x},\bar{t})=\frac{\partial L_{R}}{\partial \dot{x}}(\gm_{R}(\bar{t}),\dot{\gm}_{R}(\bar{t})).
\end{equation}
Since $R\geq \bar{R}$, then
it follows from Lemma \ref{key} that there exists $t_0>0$ independent of $R$ such that for $\bar{t}>t_0$ and any $s\in [0,\bar{t}]$,
\[H\left(\gm_{R}(s),\frac{\partial L_{R}}{\partial \dot{x}}(\gm_{R}(s),\dot{\gm}_{R}(s))\right)\leq 1.\]
Then $(\bar{x},\partial_xu_{R}(\bar{x},\bar{t}))\in \Omega$.   Moreover, for each $(x,t)\in \mathcal{E}_R$ and $t>t_0$, we have
\[|\partial_xu_{R}(x,t)|\leq \bar{R},\]
 since $\bar{R}$ is independent of  $(x,t)$. It follows that for $R>\bar{R}$, $(x,t)\in \mathcal{E}_R$ and $t>t_0$,
$u_{R}(x,t)$ satisfies
\[\partial_tu(x,t)+H_{R}(x,\partial_xu(x,t))=0.\]
 Hence, for $(x,t)\in \mathcal{E}_R$ and  $t>t_0$, $u_{R}(x,t)$ satisfies (\ref{uh}).
By Proposition \ref{semi} (a), $u_{R}(x,t)$ is a viscosity solution of (\ref{uh}). This completes the proof of Lemma \ref{tgeq}.
\End

\begin{Lemma}\label{bbddd}
Given $t>t_0$, $T^R_t\varphi(x)$ is uniformly bounded for each $R>\bar{R}$.
\end{Lemma}
\Proof
Let $\gamma_R:[0,t]\rightarrow M$ be a curve achieving the the infimum of (\ref{semmmmm}) with $\gamma_R(t)=x$. By Lemma \ref{key}, for $R>\bar{R}$, there holds
\begin{align*}
T^R_t\varphi(x)&=\varphi(\gamma_R(0))+\int_0^tL_R(\gm_R(s),\dot{\gm}_R(s))ds,\\
&=\varphi(\gamma_R(0))+\int_0^tL(\gm_R(s),\dot{\gm}_R(s))ds,
\end{align*}
which implies for any $x\in M$,
\[|T^R_t\varphi(x)|\leq \max_{x\in M}|\varphi(x)|+t\max_{(x,\dot{x})\in \Omega^*}L(x,\dot{x}).\]
This completes the proof of Lemma \ref{bbddd}.\End

By Lemma \ref{bbddd} and Lemma \ref{key}, a standard argument shows that given $t>t_0$, $T^R_t\varphi(x)$ is equi-Lipschitz for each $R>\bar{R}$ (see \cite[Proposition 5.5]{SWY}). It follows from  Lemma \ref{tgeq} that for $t>t_0$, the viscosity solution $u(x,t)$ of (\ref{hje11}) can be represented as $\liminf_{R\rightarrow\infty}T_t^{R}\varphi(x)$.
In the following, we consider the case with $t\in [0,t_0]$.
\begin{Lemma}\label{lipin}
Let $\psi(x)$ be a Lipschitz function, then there exists $\tilde{R}>0$ such that for $(x,t)\in M\times[0,t_0]$, $u_{\tilde{R}}(x,t)$ is the viscosity solution of (\ref{uh}) with $u_{\tilde{R}}(x,0)=\psi(x)$.
\end{Lemma}
\Proof
Based on uniqueness and regularity theory of viscosity solutions (\cite[Theorem 8.2]{Ba3}, \cite[Theorem 2.5]{ds}), under the assumptions (H1)-(H3), there exists a unique Lipschitz viscosity solution $u(x,t)$ of (\ref{uh}) with $u(x,0)=\psi(x)$. At the differentiable points of $u(x,t)$ on $M\times [0,t_0]$, we have
\begin{equation}
|\partial_xu(x,t)|\leq K,
\end{equation}where $K$ is a constant. Taking $\tilde{R}\geq K$, it follows from a similar argument as the one in the proof of Lemma \ref{tgeq} that for $(x,t)\in M\times[0,t_0]$, $u(x,t)$ is the viscosity solution of
\begin{equation}\label{urb}
\begin{cases}
\partial_tu(x,t)+H_{\tilde{R}}(x,\partial_xu(x,t))=0,\\
u(x,0)=\psi(x).
\end{cases}
\end{equation}
On the other hand, $u_{\tilde{R}}(x,t)$ is also a viscosity solution of (\ref{urb}). By the uniqueness of the viscosity solution of (\ref{urb}), we have $u(x,t)\equiv u_{\tilde{R}}(x,t)$ for $(x,t)\in M\times[0,t_0]$. This completes the proof of Lemma \ref{lipin}.
\End

\textbf{Proof of Theorem \ref{two}:}
First of all, we consider the case of $t\in [0,t_0]$, where $t_0$ is determined by Lemma \ref{lipin}.
For a given initial data $\varphi(x)\in C(M,\R)$, we choose a sequence of Lipschitz functions $\varphi_n(x)$ such that $\varphi_n\rightarrow\varphi(x)$ in the $C^0$-norm. Let $u_{R}^n(x,t)$ be
the viscosity solution of the following equation:
\begin{equation}
\begin{cases}
\partial_tu(x,t)+H_{R}(x,\partial_xu(x,t))=0,\\
u(x,0)=\varphi_n(x).
\end{cases}
\end{equation}
By (\ref{lo}), we have $u_{R}^n(x,t)=T_t^{R}\varphi_n(x)$. Let
\[u_n(x,t):=\liminf_{R\rightarrow\infty}T_t^{R}\varphi_n(x).\]
It follows
from Lemma \ref{lipin} that $u_n(x,t)$ is the viscosity solution of
\begin{equation}
\begin{cases}
\partial_tu(x,t)+H(x,\partial_xu(x,t))=0,\\
u(x,0)=\varphi_n(x).
\end{cases}
\end{equation}

\noindent\textbf{Claim:}
\begin{equation}\label{untt}
\lim_{n\rightarrow\infty}u_n(x,t)=\liminf_{R\rightarrow\infty}T_t^{R}\varphi(x).
\end{equation}

\noindent\textbf{Proof of the claim:}
It is easy to see that for  given $\tilde{R}>0$ and $n\in \N$,
\begin{align*}
\inf_{R>\tilde{R}}\left(T_t^{R}\varphi_n(x)-T_t^{R}\varphi(x)\right)&
\leq\inf_{R>\tilde{R}}T_t^{R}\varphi_n(x)-\inf_{R>\tilde{R}}T_t^{R}\varphi(x),\\
&\leq \sup_{R>\tilde{R}}\left(T_t^{R}\varphi_n(x)-T_t^{R}\varphi(x)\right).
\end{align*}
By virtue of the non-expansiveness of $T_t^{R}$, we have
\[\|T_t^{R}\varphi_n(x)-T_t^{R}\varphi(x)\|_\infty\leq \|\varphi_n(x)-\varphi(x)\|_\infty,\]
where $\|\cdot\|_\infty$ denotes the $C^0$-norm. Hence,
\begin{equation}
\|\inf_{R>\tilde{R}}T_t^{R}\varphi_n(x)-\inf_{R>\tilde{R}}T_t^{R}\varphi(x)\|_\infty\leq \|\varphi_n(x)-\varphi(x)\|_\infty,
\end{equation}
Since
$
\liminf_{R\rightarrow\infty}=\lim_{\tilde{R}\rightarrow\infty}\inf_{R>\tilde{R}}
$,
then we have
\begin{equation}
\|\liminf_{R\rightarrow\infty}T_t^{R}\varphi_n(x)
-\liminf_{R\rightarrow\infty}T_t^{R}\varphi(x)\|_\infty\leq \|\varphi_n(x)-\varphi(x)\|_\infty.
\end{equation}
Moreover, $u_n(x,t)$ converges to $\liminf_{R\rightarrow\infty}T_t^{R}\varphi(x)$ in the $C^0$-norm on $M\times [0,t_0]$ as $n\rightarrow\infty$, which verifies the claim (\ref{untt}).

Let $\bar{u}(x,t):=\liminf_{R\rightarrow\infty}T_t^{R}\varphi(x)$. It follows from the stability of viscosity solutions (\cite[Theorem 8.1]{F3}) that  for $(x,t)\in M\times[0,t_0]$, $\bar{u}(x,t)$ is the viscosity solution of (\ref{hje11}).

 Second, it follows from Lemma \ref{tgeq}  that for $t>t_0$, the viscosity solution $u(x,t)$ of (\ref{hje11}) can be represented as $\liminf_{R\rightarrow\infty}T_t^{R}\varphi(x)$. By virtue of the uniqueness of the viscosity solution of (\ref{hje11}) under the assumptions (H1)-(H3) \cite[Theorem 2.5]{ds}, it follows that for $(x,t)\in M\times[0,\infty)$,
\[u(x,t)=\liminf_{R\rightarrow\infty}T_t^{R}\varphi(x).\]
 In particular, there exists $t_0>0$ such that for $t>t_0$, $u(x,t)=T_t^{\hat{R}}\varphi(x)$ where $\hat{R}=\max\{\bar{R},\tilde{R}\}$. Since $T_t^{\hat{R}}\varphi(x)$ is  Lipschitz continuous and its Lipschitz constant is independent of $\varphi$ (\cite[Proposition 4.6.6]{F3}). By Lemma \ref{tgeq}, $t_0$ is also independent of $\varphi$. Hence, for $t>t_0$, $u(x,t)$ is $\iota$-Lipschitz continuous and $t_0,\iota$ are independent of $\varphi$.

So far, we have completed the proof of Theorem \ref{two}.
\End


 \vspace{2ex}
\noindent\textbf{Acknowledgement}
 The authors sincerely
thank the referees for their careful reading of the manuscript and
invaluable comments.
The authors  also would like to thank
Prof. Jun Yan for many helpful discussions.
X. Li was partially under the support of NSFC (Grant No. 11471238) and
L. Wang was partially under the support of NSFC (Grant No.  11631006, 11401107).

\vspace{2em}

{\sc Xia Li}

{\sc School of Mathematical and Physics, Suzhou University of Science and Technology,
Suzhou Jiangsu, 215009,
China.}

 {\it E-mail address:} \texttt{lixia0527@188.com}

\vspace{1em}

{\sc Lin Wang}

{\sc Yau Mathematical Sciences Center, Tsinghua University, Beijing 100084, China.}

 {\it E-mail address:} \texttt{lwang@math.tsinghua.edu.cn}

\end{document}